\documentclass[leqno
]{amsart}

\usepackage[a4paper, total={6in, 8in}]{geometry}
\usepackage[utf8]{inputenc}

\usepackage[dvipsnames]{xcolor} 
\colorlet{cite}{LimeGreen!50!Green}
\usepackage{tikz}  
\usetikzlibrary{arrows,positioning}
\tikzset{ 
  baseline=-2.3pt,
  text height=1.5ex, text depth=0.25ex,
  >=stealth,
  node distance=2cm,
  mid/.style={fill=white,inner sep=2.5pt},
}

\usepackage{lmodern}
\usepackage{tikz-cd}
\usepackage{amsmath, amsthm, amssymb, amsfonts}

\usepackage{graphicx,caption,subcaption}
\usepackage{braket}  
\usepackage{microtype}

\usepackage[%
  bookmarks=true,			
  unicode=true,			
  pdftitle={Open question about vector bundles and their moduli},		%
  pdfauthor={Elizabeth Gasparim}{Francisco Rubilar},	%
  pdfkeywords={moduli}{vector bundle},	
  colorlinks=true,		
  linkcolor=Blue,			
  citecolor=cite,		
  filecolor=magenta,		
  urlcolor=RoyalBlue			
]{hyperref}				
\usepackage{cleveref}

\newtheoremstyle{mydef}
{}		
{}		
{}		
{}		
{}	
{}		
{ }		
{\thmname{#1}\thmnumber{ #2}\thmnote{ #3}}	

\theoremstyle{plain}	
\newtheorem{theorem}{Theorem}[section] 
\newtheorem{proposition}[theorem]{Proposition}
\newtheorem*{theorem*}{Theorem}

\newtheorem{conjecture}[theorem]{Conjecture}

\theoremstyle{mydef} 
\newtheorem{definition}[theorem]{Definition}
\newtheorem*{conjecture*}{Conjecture}

\theoremstyle{remark}
\newtheorem{remark}[theorem]{Remark}
\newtheorem{notation}[theorem]{Notation}
\newtheorem{example}[theorem]{Example}

\newtheorem{question}{\bf Question}

\newtheorem*{proposition*}{Proposition}
\newtheorem*{lemma*}{Lemma}
\newtheorem*{corollary*}{Corollary}

\newcommand{\CC}{\mathbb{C}}

\newcommand{\NN}{\mathbb{N}}
\newcommand{\RR}{\mathbb{R}}

\newcommand{\ZZ}{\mathbb{Z}}

\newcommand{\PP}{\mathbb{P}}

\newcommand{\Mm}{\mathcal{M}}

\newcommand{\Oo}{\mathcal{O}}

\newcommand{\ce}{\mathrel{\mathop:}=}

\DeclareMathOperator{\Pic}{Pic}
\DeclareMathOperator{\End}{End}
\DeclareMathOperator{\Tot}{Tot}
\DeclareMathOperator{\codim}{codim}

\begin{document}

\author{E. Ballico, E. Gasparim, F. Rubilar}
\address{ E. Ballico - Dept. Mathematics, University of Trento, I-38050 Povo, Italy.
E. Gasparim {\tiny and} F. Rubilar - Depto. Matem\'aticas, Univ. Cat\'olica del Norte, Antofagasta, Chile. 
 ballico@science.unitn.it, etgasparim@gmail.com,  francisco.rubilar@ucn.cl}
\title{25 open questions about vector bundles and their moduli}

\begin{abstract}
We present  25 open questions about moduli spaces of vector bundles and related topics,
and discuss some longstanding conjectures. We hope to inspire young researchers to 
engage in this area of research.
\end{abstract}
\maketitle
\tableofcontents

We present  questions about vector bundles and their  moduli
and discuss longstanding conjectures of Iitaka, Ulrich, Mercat, Butler, and Atiyah--Jones.
This note consists of 6 independent sections:
 Iitaka dimension, Ulrich and Buchsbaum bundles,
 Higgs and coHiggs bundles, Brill--Noether theory, local numerical invariants, and 
topology of moduli spaces. It is not intended as a survey,  we favored brevity over thoroughness.
 We simply present the 
background notions necessary for stating our 25 open questions, which we believe might be useful
for young researchers interested in the study of  vector bundles and their moduli. 

The main goal of the text is to motivate research in the topic of vector bundles and their geometry,
 which has been for many years a very fruitful source of interesting and valuable mathematical  results. 
 Despite the fact that some of the conjectures discussed here have been proved in  particular cases and 
moreover some counterexamples have been found, we consider it interesting to look at the cases that remain open.
 Nevertheless, 
 since we are intending this text as a guide
  of future work for young researchers we provide many references.

\section{Iitaka dimension of vector bundles}
A vector bundle on a smooth projective variety, when it is generically generated by global sections, 
yields a rational map to a Grassmannian, called the 
Kodaira map (Def.\thinspace\ref{kod-map}). 
The Iitaka conjecture discusses the asymptotic behavior of the Kodaira maps for the symmetric powers of a vector bundle.
The theory of Iitaka fibrations and Iitaka dimension is well understood in the case of line bundles.
However, in the more general setup of higher rank bundles,
the situation becomes more involved as illustrated by 
 the following long standing fundamental question.

\begin{conjecture}[Iitaka conjecture] The Iitaka 
dimensions associated to a fibration $F\rightarrow E \rightarrow B$ 
 satisfy $\kappa(E) \geq \kappa(F)+ \kappa(B)$. 
\end{conjecture}

The conjecture was proved in dimension $\leq 6$ by Birkar \cite{ birkar} but remains open for higher dimensions
(see \cite{fujino}).
We  recall  basic definitions and  state some open questions.

If $L$ is a line bundle over a projective variety $X$, we can define a rational map
$$\phi_m\colon X\dashrightarrow Y_m\subseteq\mathbb{P}(H^0(X,L^{\otimes m})),$$
whenever $m$ is a positive integer such that $H^0(X,L^{\otimes m})\neq0$; here $Y_m$ is the closure of $\phi_m(X)$ in $\mathbb{P}(H^0(X,L^{\otimes m}))$.

Set $N(L)\ce \{m\in\mathbb{N}\colon H^0(X,L^{\otimes m})\neq0\}$. For any $m>>0$ in $N(L)$,
we have  rational maps $\phi_m\colon X\dashrightarrow Y_m$. We will see that these $\phi_m$  are birationally equivalent to the Iitaka fibration 
(Def.\thinspace\ref{def:Iitaka-fib-line-bundles}). The number $\kappa(L):=\dim Y_m$ remains constant for sufficiently large $m$ 
and is the Iitaka dimension of $L$, that interests us in this section. Formally:

\begin{definition}\label{id}
	Assume that $X$ is a normal projective variety and $L\rightarrow X$ is a line bundle. The {\bf Iitaka dimension of} $L$ is defined as
	\begin{equation}
	\kappa(L)=\kappa(X,L)\ce \max_{m\in N(L)}\{\dim\phi_m(X)\},
	\end{equation}
	provided that $N(L)\neq(0)$. If $H^0(X,L^{\otimes m})=0$ for all $m>0$, one puts $\kappa(X,L)=-\infty$. If $X$ is non-normal, 
	we pass to its normalization $\nu\colon X'\rightarrow X$ and set 
	\begin{equation}
	\kappa(X,L)=\kappa(X',\nu^*L).
	\end{equation}
	Finally, for a Cartier divisor $D$ one takes $\kappa(X,D)=\kappa(X,\Oo_X(D))$.
\end{definition}
Thus either $\kappa(X,L)=-\infty$, or else
$0\leq \kappa(X,L) \leq \dim X.$

\begin{example}[Kodaira dimension]
	Let $X$ be a smooth projective variety, and $K_X$ the canonical divisor on $X$. Then $\kappa(X)=\kappa(X,K_X)$ is the {Kodaira dimension} of $X$. The Kodaira dimension of a singular variety is defined to be the Kodaira dimension of any smooth model, that is,  a smooth 
	projective variety birrationally equivalent to $X$.
\end{example}

%
%

\begin{theorem}[Iitaka fibration]\label{it}
	Let $X$ be a normal projective variety, and $L$ a line bundle on $X$ such that $\kappa(X,L)>0$. 
	Then for all sufficiently large $m\in N(X,L)$, the rational maps
	 $\phi_m\colon X\dashrightarrow Y_m$ are birationally equivalent to a fixed algebraic fiber space 
	$$\phi_{\infty}\colon X_{\infty}\dashrightarrow Y_{\infty}$$
	of normal varieties, and the restriction of $L$ to a very general fibre of $\phi_{\infty}$ has Iitaka dimension $=0$. 
	More specifically, there exists for large $m\in N(X,L)$ a commutative diagram
	
	\begin{center}
		\begin{tikzcd}
		X  \arrow[d,"\phi_m",dashed,swap]
		& X_{\infty}\arrow[l, "u_{\infty}",swap] \arrow[d, "\phi_{\infty}"] \\
		Y_m 
		& Y_{\infty}\arrow[l,"\nu_m",dashed]
		\end{tikzcd}
	\end{center}
	of rational maps and morphisms, where the horizontal maps are birational and $u_{\infty}$ is a morphism. One has $\dim Y_{\infty}=\kappa(X,L)$. Moreover, if we set $L_{\infty}=u_{\infty}^*L$, and take $F\subseteq X_{\infty}$ to be a very general fibre of $\phi_{\infty}$, then 
	$$\kappa(F,L_{\infty}|F)=0.$$
\end{theorem}

\begin{definition}\label{def:Iitaka-fib-line-bundles}
	In the setup of  Thm.\thinspace \ref{it}, 
	$\phi_{\infty}\colon X_{\infty}\rightarrow Y_{\infty}$
	is the {\bf Iitaka fibration associated to $L$}. It is unique up to birational equivalence. 
	The Iitaka fibration of a divisor $D$ is defined by passing to $\Oo_X(D)$.
\end{definition}

\begin{definition}
	The {\bf Iitaka fibration} of an irreducible variety $X$ is by definition the Iitaka fibration associated to the canonical bundle on any non-singular model of $X$. A very general fibre $F$ of the Iitaka fibration satisfies $\kappa(F)=0$.
	
\end{definition}

Now we present an example, which is a result due to Ueno--Kawamata--K\'ollar.
\begin{example}
	Let $X\subseteq A$ be an $n$-dimensional irreducible subvariety of an Abelian variety
	and denote by $\mathbf{G}$ the Grassmannian of $n$-dimensional subspaces of the tangent space $T_0A$. 
	Then the Gauss map determines a rational map
	$$\gamma\colon X\dashrightarrow Y=\gamma(X)\subseteq \mathbf{G},$$
	which can be identified with the Iitaka fibration of $X$.
\end{example}

\begin{remark}\cite[Ex.\thinspace2.1.9]{laz}
	The Iitaka dimension is not an invariant of
	line bundles under deformations. Indeed, if $L$ is a line bundle in $\Pic^0(X)$, then $\kappa(X,L)=0$ if $L$ is trivial or torsion, but $\kappa(X,L)=-\infty$ otherwise.
\end{remark}


%
In a recent paper, Mistretta and Urbinati \cite{MU}
 generalized this construction for any vector bundle $E\rightarrow X$, and Jow \cite{Jow} gave 
 a lower bound for the Iitaka dimension of $E$, we recall the  definition.

	\begin{definition}
		Let $E$ be a vector bundle on a projective variety $X$, and let $U$ be an open subset of $X$. We say that
		\begin{enumerate}
			\item $E$ is \textbf{globally generate}d on $U$ if the evaluation map $H^0(X,E)\rightarrow E_x$ is surjective for every point $x\in U$.
			\item $E$ is \textbf{generically generated} if it is globally generated on some nonempty open subset of $X$.
			\item $E$ is \textbf{asymptotically generically generated} (AGG) if for some positive integer $m$, the $m$th symmetric power $S^mE$ of $E$ is generically generated.
		\end{enumerate}
\end{definition}

	\begin{definition}
		\label{kod-map}
		Let $E$ be an AGG vector bundle of rank $r$ on a projective variety $X$. Let $m$ be a positive integer such that $S^mE$ is globally generated on a nonempty open subset $U\subset X$.
		Set the notation $$\sigma_m(r)\ce \mbox{rank} S^mE =
		\binom{m+r-1}{m}.
		$$	
		Then we can define a rational map 
		\[\varphi_m\colon X\dashrightarrow \mathbb{G}(H^0(X,S^mE),\sigma_m(r))\]
		by sending a point $x\in U$ to the $\sigma_m(r)$-dimensional quotient $[H^0(X,S^mE)\twoheadrightarrow S^mE_x]$ of $H^0(X,S^mE)$ under the evaluation map. We call $\varphi_m$ the $m$th \textbf{Kodaira map} of $E$.
\end{definition}

Note that if $\mbox{rk}(E) =1$ then the maps $\varphi_m$ coincide with the maps $\phi_m$ used in definition \ref{id}.

\begin{theorem}\label{jow-thm-3}\cite[Thm.\thinspace3]{Jow}
	Let $E$ be an AGG vector bundle on a complex projective variety $X,$ and denote 
	$\mathbf{N}(E)=\{m\in\mathbb{N}| S^mE \textup{ is generically generated}\}.$
	For each $m\in \mathbf{N}(E)$, let $\varphi_m$ the $m$th Kodaira map of $E$, and let $Y_m$ be the closure of $\varphi_m(X)$. Then for all sufficiently large $m\in \mathbf{N}(E)$, the rational maps $\varphi_m\colon X\dashrightarrow Y_m$ are birationally equivalent to a fixed surjective morphism of projective varieties
	\[\varphi_{\mathbb{G}}\colon X_{\mathbb{G}}\rightarrow Y_{\mathbb{G}}.\]
	That is, there exists a commutative diagram
	\begin{center}
		\begin{tikzcd}
		X  \arrow[d,"\varphi_m",dashed,swap]
		& X_{\mathbb{G}}\arrow[l, "u_{\mathbb{G}}",swap] \arrow[d, "\varphi_{\mathbb{G}}"] \\
		Y_m 
		& Y_{\mathbb{G}}\arrow[l,"\nu_m",dashed]
		\end{tikzcd}
	\end{center}
	where the horizontal maps are birational and $u_{\mathbb{G}}$ is a morphism.
\end{theorem}

In the setup of Thm.\thinspace\ref{jow-thm-3}, we obtain the definition of Iitaka dimension for  a higher rank bundle.

\begin{definition}
		The \textbf{Iitaka dimension} of $E$ is
		$$\kappa(E)= \kappa(X,E)\ce \dim Y_{\mathbb{G}}.$$
	\end{definition}
	
	A very interesting result due to Mistretta and Urbinati  \cite[Rmk.\thinspace4.5]{MU} says that if $E$ is strongly semiample, then  $\kappa(X,E)=\kappa(\det E)$.

Let $X$ be a smooth projective variety. Fix numerical invariants, 
	say rank and Chern classes such that there exists a moduli space $\mathcal{M}$ of semistable vector bundles on $X$. 
\begin{question}\label{q:strat-by-iitaka-dim}
What are the properties of the 
 stratification of  the open part of $\mathcal{M}$ parametrizing stable vector bundles obtained by loci of fixed Iitaka dimension?
\end{question}

\begin{question}\label{ins}
	If $\dim X>1$, then  in the properly semistable part of $\Mm$ there are points associated to  non-locally free sheaves. 
	What is the optimal adaptation of  the  definition of Iitaka dimension to this more general situation?\end{question}

Then, once Question \ref{ins} has been solved, in the same spirit of Question \ref{q:strat-by-iitaka-dim}, we  ask:

\begin{question}\label{q:strat-of-semiestable}
What are the geometric features of the stratification  of the subset of 
$\mathcal{M}$ parametrizing equivalence classes of properly  semistable sheaves on $X$
obtained by loci of  fixed  Iitaka dimension?
\end{question}

\section{Ulrich and Buchsbaum bundles}
Let $(X,L)$ be a polarized pair, that is,  an integral projective variety $X$ (or just a projective scheme) and  an ample line bundle 
$L$ on $X$. Set $n:= \dim X$. Let $\mathcal E$ be a coherent sheaf on $X$ 
whose stalks at all $x\in X$ are $\Oo_{X,x}$-modules with depth $n$, see \cite{bes}. If $X$ is smooth, then $\mathcal E$ is locally free. 

\begin{definition} A coherent sheaf $\mathcal E$ on a polarized pair  $(X,L)$ is called 
{\bf arithmetically Cohen--Macaulay} (ACM) if the following conditions hold:
\begin{itemize}
\item[$\iota.$] $\mathcal E$ is locally Cohen--Macaulay, 
that is, the depth of the stalk $\mathcal E_x$  equals  $\dim \mathcal O_{X,x}$ for any point $x$ in $X$,
\item[$\iota\iota.$]$H^i(\mathcal E\otimes L^{\otimes t}) = 0$  for all $t \in \mathbb Z$  and $i = 1,..., \dim X -1$.
\end{itemize}
\end{definition}

Note that this definition of Arithmetically Cohen--Macaulay sheaf (or ACM) 
depends on the polarization $L$ and it is an open condition, 
hence in any moduli space $\Mm$ of vector bundles on $X$ it is
 satisfied on an open subset, which is possibly empty. Examples for surfaces appear in  
\cite{BHP} and some  examples for fourfolds are presented in \cite{casfm}. 
A stronger requirement that an ACM sheaf may satisfy is to be Ulrich. 
Ulrich bundles (when they exists) seems to be the building blocks of the derived category $D^b(X)$.
 A fundamental question concerning Ulrich bundles is stated in \ref{ulrich-on-proj-var}, but first we give a formal definition.

\begin{definition}\cite{beau1}
	Let $E$ be a vector bundle on  a smooth variety and  $X\subset\mathbb{P}^N$. We say that $E$ is an 
	{\bf Ulrich} bundle if one of the  following (equivalent) conditions is satisfied:
	\begin{enumerate}
		\item There exists a linear resolution 
		$$0\rightarrow L_c\rightarrow L_{c-1}\rightarrow\ldots \rightarrow L_0\rightarrow E\rightarrow 0$$
		with $c=\codim(X,\mathbb{P}^N)$ and $L_i=\Oo_{\mathbb{P}^N}(-i)^{b_i}$. 
		\item The cohomology $H^{\bullet}(X,E(-p))$ vanishes for $1\leq p \leq\dim(X)$.
		\item If $\pi\colon X\rightarrow \mathbb{P}^{\dim(X)}$ is a finite linear projection, then the vector bundle $\pi_*(E)$ is trivial.
	\end{enumerate} 
\end{definition}

The equivalence  of the three conditions above is  a result of Eisenbud--Schreyer--Weyman  \cite{ESW}.
 \cite{CMJ} surveys the main results and developments of the theory of Ulrich bundles in the last 30 years. Other useful references are
  \cite{ahmpl, chgs, Cs}.

%

\begin{remark}
A recent paper by F. C. Yhee \cite{Y} gives  examples of a non-normal local integral domain without Ulrich modules.
 The examples are localizations of  non-normal finitely generated integral domains over $\CC$. 
Work in progress by the first author and collaborators  is modifying the examples of Yhee to obtain
an exact characterization of integral non-normal  projective variety $X$ without Ulrich sheaves. 
\end{remark}

\begin{question}
Does there exist a  normal variety $X$ without Ulrich sheaves?
\end{question}

A very open conjecture  (on which we are skeptical but is part of the mathematical folklore)
 says:
 
 \begin{conjecture}[folklore]
   Ulrich bundles exist on any polarized pair $(X,L)$. 
 \end{conjecture}

 There are many papers describing Ulrich sheaves on specific choices of  $(X,L)$
  such as \cite{beau1,chgs,cas,casfm,lm} and references therein.
 The reader may be interested in this list, because usually such sheaves
  give good insights about the derived categories of coherent sheaves $D^b(X)$.
 Our question in this topic is:

\begin{question}
\label{ulrich-on-proj-var}	
Which smooth projective varieties carry  Ulrich bundles?
\end{question}

The original conjecture by Ulrich is the quite optimistic expectation of existence.

\begin{conjecture}[Ulrich]
Every smooth projective variety  $X\subset \PP^{N}$ has an Ulrich bundle.
\end{conjecture}

The next special types of vector bundles we  consider are 
Buchsbaum bundles.

\begin{notation}
	Let $X$ be a smooth and connected $n$-dimensional manifold with  $n\geq 2$ such that  $h^i(\Oo_X)=0$ for all $i=1,\dots ,n-1$.
	We denote by {\bf $A(X)$} the set of all ample lines bundles $L$ on $X$ satisfying $h^i(L^{\otimes t}) =0$ for all $t\in \ZZ$ and all $i=1,\dots ,n-1$.
\end{notation}	
	
	\begin{question}\label{aq1.1}
		 Under what conditions is  $A(X)$ nonempty? What is $\dim A(X)$?  
	\end{question}
	
	For a fixed  positive integer $r$ we focus now on a case when  $A(X)$ is quite large.	
	
	\begin{question}\label{aq1.2}
	 What are the properties of the subset $A(X,r)$ of all rank $r$ vector bundles $E$ on $X$ such that $h^i(E\otimes L^{\otimes t}) =0$ for all $t\in \ZZ$, all $i=1,\dots,n-1$ and all $L\in A(X)$?
	\end{question}
	
	\begin{definition}
		A vector bundle	$E$ is said to be {\bf weakly $k$-Buchsbaum with respect to a line bundle $L$},
		if for each $t\in \NN$ the evaluation map $H^0(L^{\otimes k})\otimes H^0(E\otimes L^{\otimes t})\to H^0(E\otimes L^{\otimes (t+k)})$ has maximal rank, i.e. it is either injective or surjective.	
		 Fix any subset $\Sigma \subseteq A(X)$. $E$ is called
		  {\bf weakly $k$-Buchsbaum with respect to $\Sigma$} if it is weakly $k$-Buchsbaum for all $L\in \Sigma$. 
	\end{definition}
	
	For the varieties  $X$  considered here, duality gives $h^i(\omega _X)=0$ for all $i=1,\dots,n-1$. 
	If $X$ is a product of projective spaces, then by Kunneth's formula, $A(X)$
	equals  the set of all ample line bundles.  As consequence of Horrock's theorem, over a projective space   $\PP^n$ with $n\ge
	2$, we have that $E\in A(\PP^n,r)$ if and only if it splits, that is,  $E\cong \oplus _{i=1}^{r}\Oo_{\PP^n}(a_i)$ for some $a_1\ge \cdots \ge a_r$.

 There is at the moment no complete classification of $k$-Buchsbaum  vector bundles, not 
	even if we restrict to   complex manifolds $X$ of an elementary nature.	So, a natural question is:
	
	\begin{question}\label{aq2}
		What properties characterize the set of all weakly $k$-Buchsbaum vector bundles on a $n$-dimensional complex manifold $X$?
	\end{question}
	
	We suggest starting with the cases of $\PP^a\times \PP^{n-a}$ and $(\PP^1)^n$ and comparing them.

\section{Higgs and coHiggs bundles}

There exist many concepts of decorated bundles. One of the motivations to adding decorations is the search for  
fine moduli spaces. Classical concepts of decoration are stable, semistable, and parabolic bundles, which belong to a vast literature 
\cite{Ha1,Ha2,HL, MS, NR, NS}. 
More recent decorations were obtained by adding Higgs bundles,
with the name being motivated by the physics concepts of a Higgs field. 
Higgs bundles were introduced by N. Hitchin \cite{Hi1,Hi2,Hi3,Hi4} and a huge number of papers  explores them. The interested reader may find a partial list in the recent \cite{bll}. The case of coHiggs sheaves has been  less studied, because they are much less common, as  Proposition \ref{less} shows. A good reference for the state of art on
 both Higgs and coHiggs bundles is \cite{Swo}.
We first recall basic definitions and then state open questions.


\begin{definition}\cite{GGB} Let $X$ be a smooth complex projective variety.
\begin{itemize}
\item	A {\bf Higgs bundle} on $X$ is a pair $(E,\varphi)$, where $E$ is a holomorphic vector bundle and $\varphi\colon E \to E\otimes \Omega ^1_X$, the Higgs field, is a holomorphic $1$-form with values in $\End(E)$ such that $\varphi\wedge \varphi=0$.
\item	A {\bf coHiggs bundle} on $X$ is a pair $(E,\varphi)$, where $E$ is a holomorphic vector bundle and $\varphi\colon E\to E\otimes T_X$ the coHiggs field, is a holomorphic vector field with values in $\End(E)$ such that $\varphi\wedge \varphi=0$.
\end{itemize}	
\end{definition}

\begin{remark}
When $\dim X=1$ the integrability condition is automatically satisfied.
\end{remark}

\begin{example}
	In the simplest examples $E$ is a line bundle and $\varphi$ is a holomorphic $1$-form.
\end{example}

\begin{remark}
	One of the most important features of Higgs bundles is that they have continuous moduli, i.e., they come in families parametrized by the points of a geometric space (in fact, a quasi-projective variety) known as their moduli space.
\end{remark}

\begin{example}
Moduli spaces of Higgs bundles on Riemann surfaces are noncompact hyperk\"ahler manifolds. 
\end{example}

In the subject of Higgs bundles and coHiggs bundles we emphasize
 the so called logarithmic case, i.e. the case of a pair $(X,D)$ with $D$ effective divisor, e.g. with simple normal crossing.
When $X$ is a smooth projective curve and $D$ is a finite set, the setup $(X,D)$ is called a punctured Riemann Surface.
In such a case, without assuming the integrability condition, we get the definitions of Higgs and coHiggs fields. We always assume $f\ne 0$. A Higgs (or coHiggs) subsheaf is a Higgs (resp. coHiggs) sheaf $(G,h)$ with $G\subset E$ and $h = f_{|G}$. Any polarization on $X$ gives a notion of stable and semistable Higgs (or coHiggs) sheaf for which there are moduli spaces. 

\begin{proposition}\label{less}\cite[\S\thinspace 3]{R2}
Semistable bundles with nontrivial coHiggs structure do not exist on smooth curves of genus $\geq2$.
\end{proposition}

Roughly speaking, the existence of stable coHiggs bundles determines the position of $X$ towards
 a negative direction in the Kodaira spectrum. For instance, M. Corr\^{e}a proved that if $\dim(X) = 2$, the existence of a semistable coHiggs bundle $(E,f)$ of rank $2$
  with $f$ nilpotent, implies that $X$ is either uniruled, a torus, 
  or some specific elliptic surface, up to finite \'{e}tale cover \cite{Correa}. 
  Nevertheless, they are interesting varieties as attaining some of the extremal cases for the generalized geometric structures considered by N. Hitchin and M. Gualtieri \cite{Gual,Hi3, Hi4}. The following  is a partial list of the papers studying them and some of extensions of this notion, 
 such as logarithmic coHiggs sheaves \cite{VC1,BH1,BH2,BH3,BH4, BHM, Correa,Gual,Hi2,R1,R2,Rayan}.
 
There is an  
obvious extension of the notion of Higgs or coHiggs field to a collection of fields; namely,
 fixing an integer $r\ge 2$ and take $r$ Higgs or coHiggs fields $f_1,\dots ,f_r$ on the same manifold,
which then leads associated notions of semistability and stability and moduli spaces. 
Requiring  additionally that $f_1,\dots ,f_r$ be linearly independent 
we obtain  equivalent notions. 
Alternatively,  taking  $r$-dimensional linear subspaces of 
$H^0(End(E)\otimes \Omega ^1_X)$ or $H^0(End(E)\otimes T_X)$ also leads to  equivalent definitions.
However  a requirement of integrability ought to be added.
In case $\dim X>1$ various inequivalent concepts of integrability are possible.
 We  prefer the strongest one: impose that $f_i\wedge f_j=0$ for all $i,j$. 
 For instance, this is  the case considered for the 2-nilpotent coHiggs structures studied 
 in \cite{BH2,Correa} and references therein.
The next step is to upgrade from fields to bundles, hence the following question.

\begin{question}
	What are the optimal  extensions of  the concepts of Higgs and coHiggs bundles to include $r$ fields in place of $1$ in
	 the case of manifolds of dimension greater than one?
\end{question}

\begin{remark}
If $\dim X>1$ it is much easier to construct Higgs or coHiggs fields than Higgs or coHiggs bundles. For a fixed vector
bundle $E$ on $X$ the possible Higgs or coHiggs fields form the vector space $H^0(Hom (E,K_X\otimes E))$ or $H^0(E,T_X\otimes
E)$, while the integrable ones are cut out by nonlinear equations. In the examples, say $X=\PP^2$, $E$ split or $T\PP^2$,
these additional equations are quadratic ones. Thus nonexistence theorems are harder to prove for Higgs or coHiggs fields. 
Some nonexistence results follow directly from 
\cite[Thm\thinspace 1.1]{BH1}, \cite[Prop.\thinspace 2.8]{BH1}, and \cite[\S 3]{BH2}, where they are stated together with 
considerations about the implications to the nonexistence  of coHiggs fields. But certainly
the interested reader may find many other nonexistence theorems for coHiggs fields by considering other rational
complex projective manifolds. We only pose a nonvague question in
the case of multiprojective spaces $\PP^{n_1}\times \cdots \PP^{n_k}$,
 in which case the main difficulty is given by the  Picard group being isomorphic to $\mathbb {Z}^k$, 
 but broader situations are also of interest. 
\end{remark}

Let $X =\PP^{n_1}\times \cdots \times \PP^{n_k}$, $k\ge 2$, be a multiprojective space. A natural question is:

\begin{question}\label{pol}
How is the moduli space of low rank coHiggs bundles on $X$ affected by a change of polarization?
\end{question}

Recall that the Segre--Veronese embedding of $X$ associated to $\Oo_X(a_1,\dots ,a_k)$
 is related to partially symmetric tensors. From the point of view  of tensors, the most important polarization is the one with $a_i=1$ for all $i$.
 For geometry, the most important polarization is the anticanonical one (with $a_i=n_i+1$ for all $i$). 
 Thus, if $n_i\ne n_j$ for some $i, j$, 
 there is no preferred choice of polarization, hence a solution to question \ref{pol} is fundamental.

\section{Brill--Noether theory}\label{Sc}

We consider the Brill--Noether theory of stable and semistable 
vector bundles on a smooth and connected projective curve. 
The case of singular, but integral curves, mainly  nodal curves, is discussed in work of 
U. N. Bhosle \cite{bhos2,bhos3} and references therein.
The case of reducible curves appeared  in \cite{bhos1,bf1,bf2} and we will discuss it 
 in subsection \ref{Se}. Further classical results of Brill--Noether theory appear in 
\cite{ln2,ln3,ln4,lnp,lns}.

Assume  $X$ is a curve of genus $g\ge 2$ and let $\Mm(X,r,d)$ 
denote the moduli space of stable rank $r$ vector bundles on $X$ with degree $d$.
 The variety $\Mm(X,r,d)$ is smooth and irreducible of dimension $r^2(g-q)+1$. 
 
 \begin{definition}The rank $r$ {\bf Brill--Noether locus} is defined as
  $$W^{k-1}_{r,d}(X)\ce \{E\in \Mm(X,r,d):  h^0(E)\ge k \}$$ 
and  the {\bf Brill--Noether number} is defined as
  $$\rho^{k-1}_{r,d}\ce r^2(g-1) -k(k-d+r(g-1)).$$  
 \end{definition}
 
 If rank $r=1$ and $\rho ^{k-1}_{1,d}\ge 0$, then
 $ W^{k-1}_{1,d}(X)\ne \emptyset$. Fundamental results of  Brill--Noether--Gieseker--Petri 
 describe the case when $r=1$   \cite[Ch.\thinspace VIII]{acgh} showing that for a general variety $X$:
 \begin{itemize}
 \item   $W^{k-1}_{1,d}(X)\ne \emptyset$ if and only if $\rho ^{k-1}_{1,d}\ge 0$,
 \item   $\dim W^{k-1}_{1,d}(X)= \rho ^{k-1}_{1,d}$, and
 \item  if $\rho ^{k-1}_{1,d}>0$ then  $W^{k-1}_{1,d}(X)$ is irreducible.  
\end{itemize}

For $r\ge 2$ the situation is more complicated even for a general $X$, and  
 25 years ago P. E. Newstead asked for a complete description of these loci
  (nonemptiness, dimension, irreducibility when of positive dimension, and smoothness as in the Gieseker--Petri case, i.e. $W^k_{r,d}(X)$ is the singular locus of $W^{k-1}_{r,d}$). These expected statements were not quite true  (even the existence/nonexistence part).
Many results if this direction are known, and the list of references is long, with many results obtained for very low genera or for low rank \cite{bgn,hhn,ln2,ln3,ln4,lnp,lns,M2}. 

\begin{question} What is the relation between the Brill--Noether theory for line bundles on $X$ 
and the Brill--Noether theory for stable rank $r\ge 2$ vector bundles on $X$?
\end{question}

We now discuss the concept of Clifford indices, and then propose a comparison with Brill--Noether numbers.
The Clifford index $\gamma_C$ (or Cliff($C$)) of a curve $C$ measures the complexity of the curve inside its moduli
space. In \cite{ln1} Lange and Newstead proposed the following definition as a generalization of $\gamma_C$ for
higher rank vector bundles. 

\begin{definition}
Let $ E \in $  be a semistable vector bundle of rank $n$ and degree $d$
on a curve $C$ of genus $g$. The \textbf{Clifford index} of $E$ is defined as :

$$\gamma(E) \ce \mu(E) -\frac{2}{n} h^0(C, E) + 2 \geq 0$$
where $\mu(E)=d/n$, and the \textbf{higher Clifford indices} of $C$ are defined as:
$$\mbox{Cliff}_n(C) = \min\{\gamma(E) :  E \in  \mathcal U_C(n,d),  d < n(g -1), h^0(C, E) \geq 2n\}.$$

\end{definition}

We recall an old conjecture by V. Mercat concerning Clifford indices,    originally stated in \cite{M2}.

\begin{conjecture}[Mercat]
Let $E$ be a semistable vector bundle of rank $n$ and degree $d$.
\begin{enumerate}
	\item If $\gamma_1+2\leq  \mu(E)\leq 2g-4-\gamma_1$, then $h^0(E)\leq \frac{d-\gamma_1 n}{2}+n$.
	\item If $1\leq \mu(E)\leq \gamma_1+2$, then $h^0(E)\leq \frac{1}{\gamma_1+1}(d-n)+n$.
\end{enumerate}
\end{conjecture}

Even though many counterexamples are known, see for instance
  \cite{Se}, the conjecture has been proved in some cases,  \cite{BF}, so still it remained unclear whether a weaker conjecture may describe well 
  the situation of high genus. 
After the first counterexamples, the conjecture was solved  in some low rank cases for
 curves with very low gonality (i.e. the lowest degree of a nonconstant rational map to $\mathbb P^1$) with respect to their genus. 
 Subsequently, this line of  study was interrupted, while
complete descriptions  for cases of low genera appeared in papers of Lange and Newstead.
 It would be interesting to understand the relation between the 2 lines of work. The question we propose is:

\begin{question} 
Fix rank $r=2$ and high genus. What is the relation between the Brill--Noether loci and the loci of fixed Clifford indices?
\end{question}

\subsection{Coherent systems and triples}

The concept of coherent system was developed as a very strong tool 
for the study of the Brill--Noether theory of vector bundles and related topics
such as  Butler's conjecture (considered in  subsection \ref{Sd}). Coherent system is  an interesting 
concept even for $g=0,1$. Let $X$ be a smooth and connected projective curve of genus $g$. 

\begin{definition}\cite{bgmn,Nitsure}
A {\bf coherent system} of multidegree $(r,d,v)$ on $X$ is a pair $T=(E,V)$, where $E$ is a rank $r$ 
vector bundle on $X$ of degree $d$ 
and $V$ is a $v$-dimensional linear subspace of $H^0(E)$. 
A {\bf coherent subsystem} $(E',V')$ of $E$ is a coherent system $(E',V')$ with $E'$ a subsheaf of $E$ and $V'\subseteq V\cap H^0(E)$. 
Note that we do not assume that $E'$ is a subbundle of $E$. 

For each real number $\alpha \ge  0$ the {\bf $\alpha$-slope} $\mu _\alpha(T)$ of $T$ is the real number
$\frac{d}{r} +\alpha\frac{v}{r}$. 
\end{definition}

Using the $\alpha$-slope of $T$ and of its coherent subsystems we obtain
concepts  of $\alpha$-stability and $\alpha$-semistability for coherent systems. 
Let $U(X,r,d,\alpha)$ denote the moduli space of (equivalence classes of) $\alpha$-semistable coherent systems of type $(r,v)$. For fixed $r$, $d$ and $v$ there are finitely many positive real numbers $\alpha_1<\cdots <\alpha _t$ (called the critical values) such that proper $\alpha$-semistability may occur only if $\alpha \in \{\alpha_1,\dots ,\alpha_t\}$. 
In such a case $U(X,r,d,\alpha)$ 
is constant when $\alpha <\alpha_1$,  in the open intervals $(\alpha _i,\alpha _{i+1})$,
and  when $\alpha > \alpha _t$. Therefore, 
$$\lim _{\alpha \to 0^+} U(X,r,d,\alpha) \quad \text{and} \quad 
\lim _{\alpha \to +\infty} T(X,r_1,r_2,d_1,d_2,\alpha)$$
 are well-defined and very interesting. When $\alpha$ crosses $\alpha_i$ the moduli space may change, but usually in a controlled way.

A {\bf holomorphic triple} on $X$ is a triple $T=(E_1,E_2,f)$, where $E_1$ and $E_2$ are vector bundles on $X$ and $f: E_2\to E_1$, see \cite{bg,bgg,PP,PP0, BR}. 
A subtriple of $T$ is a triple $(E'_1,E'_2,f')$, where $E'_i$ is a subsheaf of $E_i$, $f(E'_2)\subseteq E'_1$ and $f' =f_{|E'_2}$.
For each real number $\alpha >0$ the $\alpha$-slope $\mu_\alpha(T)$ of $T$ is the real number
$$\mu_\alpha(T):= \frac{\deg(E_1)+\deg(E_2)+\alpha\mbox{rk}(E_2)}{\mbox{rk}(E_1)+\mbox{rk}(E_2)}.$$

Using the $\alpha$-slopes of all subtriples of $T$ we get the notions of  $\alpha$-stable and $\alpha$-semistable triples.

 Let $T(X,r_1,r_2,d_1,d_2,\alpha)$ denote the set of (equivalence classes of) holomorphic triples with ranks $r_1,r_2$ and degrees $d_1,d_2$. As in the case of coherent systems  for fixed $d_1$, $d_2$, $r_1$ and $r_2$ there are finitely many positive rational numbers $\alpha_1<\cdots <\alpha _t$ (called the critical values) such that proper $\alpha$-semistability may occur only if $\alpha \in \{\alpha_1,\dots ,\alpha_t\}$. 
Similarly,  $T(X,r_1,r_2,d_1,d_2,\alpha)$ is constant when $\alpha <\alpha_1$,  in the open intervals $(\alpha _i,\alpha _{i+1})$,
 and when $\alpha > \alpha _t$. Therefore,  
 
  $$\lim _{\alpha \to 0^+} T(X,r_1,r_2,d_1,d_2,\alpha) \quad \text{and} \quad \lim _{\alpha \to +\infty} T(X,r_1,r_2,d_1,d_2,\alpha)$$
   are well-defined and very interesting. The rational numbers $\alpha _i$ are easily computed, because they are the only positive real numbers for which a properly semistable triple with numbers $r_1,r_2,d_1,d_2$ may have a proper subsheaf with the same slope.
Obtaining  all $\alpha_i$ is done in steps,  first checking all $(k_1,k_2)$ with $k_i\in \{0,1,\dots ,r_i\}$, $(k_1,k_2)\notin \{(0,0),(r_1,r_2)\}$ 
and so on.
When $\alpha$ crosses $\alpha_i$ the moduli space may change, but usually in a controlled way. 

A detailed   study of the sets $S(\alpha_i)$ of the equivalence classes of properly $\alpha_i$-semistable triples
is needed. These are the walls that separate the moduli space of $\alpha_i$-semistable bundles. Crossing $\alpha _i$ on the left or on the right, the open chambers are again $\alpha$-stable, while some of the walls disappear (they go to $\alpha$-unstable bundles) and some of the walls are promoted to $\alpha$-stable triples, so that two chambers may merge.
The cases $g=0,1$ are studied in \cite{PP0,PP}, we ask for the details in  the next case.
\begin{question}
What is the configuration of  walls and what are the limits to $+\infty$ and to $0^+$ for genus $2$ curves?  
(it should not  depend on the genus $2$ curve).
\end{question}

  In general one may associate a moduli space to any quiver of maps between vector bundles on the same scheme, 
e.g.  \cite{ACG}.
Coherent systems and holomorphic triples were considered in the setup of logarithmic coHiggs in \cite[\S 7]{BH3}.

\subsection{Butler's conjecture}\label{Sd}
Let $X$ be a projective scheme and $E$ a rank $r$ vector bundle on $X$. Set $n:=\dim X$. 
Assume that $h^0(E)\ge n+r$ and $E$ is globally generated.
 For each integer $v$ such that $v\ge n+r$ let $G(E,v)$ denote the set of all $v$-dimensional linear subspaces of $H^0(E)$.
 An easy lemma ascribed to Serre says that a general $W\in G(E,v)$ spans $E$, i.e.
  the evaluation map $e_{E,W}\colon  W\otimes \Oo_X\to E$ is surjective.
   Thus $M_{E,W}\ce \mathrm{ker}(e_{E,W})$ is a rank $v-r$ vector bundle on $X$. Now assume $\dim X=1$. More than 25 years ago D. C. Butler conjectured (Conj.\thinspace \ref{butler-conj}) that $M_{E,W}$ is semistable
    when  $X$ is general and $E$ is general in $W^{v-1}_{r,d}(X)$ \cite{butler}. 
 Many particular cases were solved, but not all, and no counterexample was found in the smooth case, see \cite{bmno} and references therein.
 In \cite[Thm.\thinspace 5.1]{BBN} conjecture \ref{butler-conj} was proven for
line bundles on a general smooth curve. 
 Nevertheless, for a   smooth projective curve $X$ of genus $g\geq 2$ and $E$ a rank $r>1$ vector bundle on $X$, the 
 following conjecture remains open. 
 
 \begin{conjecture}[Butler]\label{butler-conj}
For $C$ a general curve of genus $g\geq 1$ and a general choice of $(E,V)$ the bundle $M_{E,V}$ is semistable.
\end{conjecture}

%
%

\subsection{Reducible curves}\label{Se}
Now we consider   a reduced and connected  projective curve $X$ with $s\ge 2$ irreducible components $X_1,\dots ,X_s$.
 For each $S\subseteq \{1,\dots ,s\}$ set $X_S\ce \cup _{i\in S} X_i$ and
$X_{{S}^c}\ce \cup_{i\notin S} X_i$. Thus $X_\emptyset =X_{{\{1,\dots ,s\}^c}} =\emptyset$, 
while in all other cases the connectedness of $X$ implies that the scheme 
$Z_S\ce X_S\cap X_{{S}^c}$ is a nonempty zero-dimensional scheme. 
We will always regard $Z_S$ as a scheme, not just a set, 
but if $X$ is nodal then each $Z_S$ has a reduced structure. 
Let $j_S\colon X_S\to X$ denote the inclusion. Set $z_S\ce \deg (Z_S)$. Obviously $Z_S =Z_{{S}^c}$.

\begin{definition}
 The {\bf connectedness degree} of $X$ is the minimum of all integers $z_S$, with $S\ne \emptyset$, $S\ne \{1,\dots ,s\}$. 
 \end{definition}

In the case of other singular varieties,  to obtain  the strongest results
it is  neecessary to allow  non-locally free coherent sheaves, as usual. 
Let $F$ be a coherent sheaf on $X$. $F$ is said to have depth $1$ if all of its nonzero stalks  have depth 1. 
Obviously a vector bundle has depth $1$. Any nonzero subsheaf of a depth $1$ sheaf has depth $1$.
To speak about stability we need a polarization on $X$, i.e. a vector $w=(w_1,\dots ,w_s)\in \RR^s$ with $w_i>0$ for all $i$. 
To obtain good moduli it is not sufficient to fix the rank of the vector bundle $E$ and the multidegree $(d_1,\dots ,d_s)$ of $E$, 
it is necessary to use all $z_S$ and their relation with $w_1,\dots ,w_n$ and the multidegree of $E$. We suggest \cite{bfmv} as an inspiration.
 To start with the concepts of  gonality and Clifford index for a stable reducible curve need to be clarified. 
 
\begin{remark}[Counterexample of Brivio--Favale to Butler's conjecture]
Let $X$ be a curve with two smooth components and a single node, i.e. assume $X$ nodal, $s=2$ and $z_1=1$. 
Thus, the arithmetic genera satisfy $p_a(X) =p_a(X_1)+p_a(X_2)$ and the integers $(p_a(X_1),p_a(X_2))$ are the only numerical invariants of $X$. If $p_a(X_i)\ge 3$ for some $i$, even the Brill--Noether theory of line bundles on $X$, may depend  on the isomorphic classes of $X$. 
S. Brivio and F. Favale \cite{bf2} proved that for this reducible curve Butler's conjecture fails quite often, in non-pathological cases.\end{remark}

\begin{question}
What assumptions in the case of a reducible curve imply a positive answer to Butler's question?
\end{question}

In particular, for the  counterexample of Brivio--Favale,
we might need to assume  that the rank is small with respect to the integer $z_1$ and choose a suitable polarization $(w_1,w_2)$.

\section{Local numerical invariants}

Consider the situation of a vector bundle $E$ (or a sheaf) on a variety $X$ and 
a subvariety $Y \subset X$. We wish to study local numerical invariants of $E$ around $Y$
such as local characteristic classes, that depend only on the behaviour of $E$ on a small 
(analytic) neighborhood $N$ of $Y$. The neighborhood can be intuitively 
though of as the normal bundle of $Y$ inside $X$. 

We start by stating  an open question of very large scope:

\begin{question}
What should be  the axioms for  a general theory of local characteristic classes of vector bundles around a contractible subvariety?
\end{question}

Once  a specific concept of local invariant has been chosen,   three other questions follow immediately. 
\begin{question}
What are its admissible numerical values?
\end{question}

That is, one wants to have bounds for the invariant in terms of the restriction of $E$ to $Y$. Then, we also 
want to know which intermediate values actually occur for some bundle, giving a question about existence.

\begin{question}
Which admissible values of the local invariant are realized by some vector bundle?
\end{question}
 
 Once existence is obtained for a fixed value, there comes the corresponding question of moduli.
 \begin{question}
 What is the moduli space of all vector bundles on $N$ with  prescribed local numerical invariants?
 \end{question}
 
 \begin{remark} Questions about local numerical invariants have many applications to 
 mathematical physics. See for instance applications to the   local  charges of instanton  
presented in \cite{GKM,GO, GSu}.
 \end{remark}

One of the simplest concepts of local invariant is the  local holomorphic Euler characteristic
around a contractible subvariety. We recall the definition.

\begin{definition}\label{localEuler}\cite[Def.\thinspace 3.9]{Bl} 
 Let $\pi\colon (Z,\ell) \rightarrow (X,x) $ 
  be a resolution of an isolated quotient singularity, $\mathcal F$ a reflexive sheaf on $Z$ and $n \ce \dim X$. 
   The local holomorphic Euler characteristic of $\pi_*\mathcal  F$ at $x$ is 
$$\chi(x, \pi_*\mathcal F) \ce \chi(\ell, \mathcal F) \ce h^0(X, (\pi_*\mathcal F)^{\vee\vee}/\pi_\mathcal *F) + 
\sum_{i=1}^{n-1} h^0(X,R^i\pi_*\mathcal F) .
$$

\end{definition}

When $X$ is a compact orbifold, Blache \cite{Bl} showed that the global Euler characteristics of $X$
 and its resolution are related by

$$\chi(X, (\pi_*\mathcal F)^{\vee\vee}) = \chi(Z, \mathcal F) + \sum_{x\in Sing \, X} \chi(x, \pi_*\mathcal F).$$

The local holomorphic Euler characteristic has been intensively studied in 
the case of a local surface containing a contractible rational surface. For  the surfaces 
$$Z_k\ce \Tot(\mathcal O_{\mathbb P^1}(-k))$$ 
local $\chi$ was studied  in 
the context of computational algebraic geometry \cite{GSw,BGK2},  singular varieties \cite{GK},
 applied to questions of the  moduli spaces of vector bundles  \cite{Ga,BGK1,BeG, BrG}, to 
 bundles on blow-ups \cite{GA1,BG5} and  to 
 questions of existence and decay of instantons in mathematical physics \cite{GKM,GO, GSu}.
 Nevertheless, even in such case, fundamental questions of existence remain open.
 
 For the case of $Z_1$, that is, the blow-up of a smooth point, \cite{BG1, BG2} gave sharp bounds and 
 proved existence 
 of vector bundles with all admissible values of numerical invariants. 
 Questions of existence are essential to the study of moduli. For instance,  for the rank 2 case, \cite[Thm.\thinspace 4.15]{BGK2} 
 showed that the pair of invariants $({\bf w}, {\bf h})\ce ( h^0(X, (\pi_*\mathcal F)^{\vee\vee}/\pi_\mathcal *F), h^0(X,R^1\pi_*\mathcal F)) $
stratifies the moduli stacks of bundles on $Z_k$ into Hausdorff components, whereas their sum, the local $\chi$, 
does not by itself produce such a stratification. 
Motivated by their applications to the  physics of instantons, the local invariants \textbf{w} and \textbf{h} are called width and height respectively.
The following question has only been completely solved in the case of  $k=1$.
\begin{question}\label{inv} Let $({w}, { h})$ be  admissible values of  local invariants  on $Z_k$, 
that is, they satisfy the bounds of  \cite[Cor.\thinspace 2.18]{BGK2}. Does there exist a vector bundle $E$ on 
$Z_k$ with $({\bf w}(E), {\bf h}(E))=({ w}, { h})$?
\end{question}

A solution to question \ref{inv} will have applications to question of 
existence and decay of instantons via the Kobayashi--Hitchin correspondence. 
Other open questions about these local numerical invariants in the context of singular varieties are 
stated in  \cite{GK}.

The next cases to study are those of local threefolds. There are two clear directions 
of study. That of the total space of a line bundle on a surface and that of a rank-2
bundle on a curve. The former was considered  in \cite{BG3, BG4} but  otherwise has not 
been sufficiently explored, whereas the latter has been intensively studied, especially in the case of 
a 3-dimensional neighborhood of a contractible  curve. In particular, a type of duality occurs 
between the 2 and 3 dimensional cases, providing isomorphisms of moduli spaces of vector bundles,
 see \cite{ABCG}.

One of the most important cases for us is that of Calabi--Yau threefolds, 
which in the case of a contractible curve leads us to the following situation. 
Contraction of a rational curve $\ell$ on a  Calabi--Yau threefold $W$ 
may occur in 3 cases \cite{Ji}, namely $N_{\ell/W}$ must be isomorphic to
one of
\begin{equation}
W_1\ce \mathcal{O}_{\mathbb{P}^1}(-1) \oplus \mathcal{O}_{\mathbb{P}^1}(-1)
\text{ ,}\quad W_2\ce \mathcal{O}_{\mathbb{P}^1}(-2) \oplus \mathcal{O}_{\mathbb{P}^1}(0)
\text{ , \ or \ } W_3\ce \mathcal{O}_{\mathbb{P}^1}(-3) \oplus \mathcal{O}_{\mathbb{P}^1}(+1)
\text{ .}
\label{3folds}
\end{equation}

Let us denote by $W_k$ any of the previous Calabi--Yau threefolds in (\ref{3folds}). 
We present some open questions whose solutions would lead to a better understanding of local numerical invariants for threefolds.

\begin{question}\label{q:lloc-inv-wk}
What  local invariants   stratify the moduli 
of holomorphic vector bundles on $W_k$ for $k\geq1$ into Hausdorff components (in the analytic topology)?
\end{question}

Once question \ref{q:lloc-inv-wk} is resolved, then we have:
\begin{question}
	Construct a Hausdorff stratification of the moduli stack $\mathfrak{M}_n(W_k)$ of holomorphic bundles on $W_k$ with $\chi(\ell,E) = n$.
\end{question}

While this question is partially solved for $W_1$, the situation for $k\geq2$ is widely open.  \cite[Thm.\thinspace5.1]{GK}
 gives sharp bounds for the local holomorphic  Euler characteristic $\chi(\ell,\pi_*E)$ in the case of rank $2$ vector bundles on $W_1$,
and  \cite[Thm.\thinspace5.3]{GK} gives bounds for $h^1(W_1,End(E))$. A natural question is:

\begin{question}
Is every admissible value of $h^1(W_1,\mathcal{E}nd(E))$ realized by some vector bundle $E$?
\end{question}

We stop here, having given a few basic questions. Nevertheless, an entire theory of local characteristic classes ought to be developed. 
We hope that these questions motivates some readers towards such a task.

\section{Topology of moduli spaces}
The study of 
moduli of vector bundles 
is a central topic in 
algebraic geometry and related areas such as algebraic topology 
and mathematical physics. 
A fruitful application to physics is the study of instanton moduli spaces, which 
connects to the study of moduli of bundles via the 
celebrated
Kobayashi--Hitchin correspondence, summarized as:
$$\left\{
\begin{array}{c}
\textup{irreducible } \mathrm{SU}(2)-\textup{instantons}\\
\textup{instantons of charge } n
\end{array}
\right\}\leftrightarrow\left\{
\begin{array}{c}
\textup{stable }\mathrm{SL}(2,\mathbb C)-\textup{bundles}\\
\textup{with } c_2=n
\end{array}
\right\}.
$$

The literature on instanton moduli contains strong results,
such as the seminal works of Donaldson \cite{D1, D2}. 
Atiyah and Jones  \cite{AJ} proved a fundamental result 
about the homology of moduli spaces of $\mathrm{SU}(2)$-instantons on the sphere 
$S^4$, namely they showed that the inclusion of the moduli space of instantons 
of charge $k$ into the moduli space of all connections modulo gauge equivalence 
is a surjection  in degrees lower than $k$. The result promotes the 
spirit of Morse theory to infinite dimensional spaces, showing that  homology in low degree is detected by 
the critical sets of the Yang--Mills functional. 
\cite{AJ}    conjectured  that this inclusion induces isomorphisms  in
homology and homotopy for the same low degree range. 
The statement is known as the Atiyah--Jones conjecture  
and the question has been generalized  to any $4$-manifold $M$. 
Let us denote by $\mathcal{M}_k$ the moduli space of $SU(2)$ instantons on $X$ of charge $k$. 
For a  principal bundle $P_k$ with $c_2=k$  on $X$ let  $\mathcal{B}_k$ be the space of gauge-equivalence classes of connections on $P_k$. Then we can state a generalized version of the Atiyah--Jones conjecture as follows.

\begin{conjecture}[Generalized Atiyah--Jones conjecture]\label{gen-AJ-conj}
Let $X$ be a $4$-dimensional real manifold. Then the inclusion $\theta_k\colon \mathcal{M}_k\rightarrow\mathcal{B}_k$ induces isomorphisms in homology and homotopy, that is,
\begin{itemize}
	\item $(\theta_k)_t\colon H_t(\mathcal{M}_k)\rightarrow H_t(\mathcal{B}_k),$ and
	\item $(\theta_k)_t\colon \pi_t(\mathcal{M}_k)\rightarrow \pi_t(\mathcal{B}_k),$
\end{itemize}
are isomorphisms for $t\leq q=q(k)$.
\end{conjecture}

The Atiyah--Jones conjecture was proven in some cases:
for $\mathrm{SU}(2)$ instantons on $S^4$   by Boyer, Hurtubise, Mann, and Milgram \cite{BHMM},
for $\mathrm{SU}(n)$ instantons on $S^4$ by Tian \cite{T}, 
for ruled surfaces by Hurtubise and Milgram \cite{HM}, and  for rational surfaces  by Gasparim \cite{Ga}. 
The conjecture remains open in all other cases. In particular, the case of irrational surfaces is missing.

\begin{question}
\label{q:general-AJ}
Is the Atiyah--Jones conjecture true for any irrational complex surface?
\end{question}

The conjecture  compares the topology of a moduli space of vector bundles with fixed Chern classes with the 
topology of the space obtained in the limit when   the top Chern class goes to infinity. 
Such  studies in general require  a choice of stability and an analysis of the subset of unstable bundles.
Our results in \cite{BGRS} argue that the intermediate step of choosing stability is negligible for the understanding of the stable limit,
simplifying the technical steps  to compare the topology of the moduli spaces to that of the corresponding stacks. 
(For the basic theory of stacks see for instance \cite{Go,O1,O2} and for  moduli of bundles on stacks  see \cite{nir, su}.) 

In \cite{BGRS}, we compare the homology of the moduli space $\mathfrak{M}$ 
of rank $2$ vector bundles  on certain surfaces of general of general type to 
the irreducible component  the moduli stack $\mathfrak{M}^s$  of vector bundles 
which contains stable  bundles. We prove the following result.

\begin{theorem}\cite[Thm.\thinspace 18]{BGRS}\label{lim}
Assume $c_2>>0$, and let $\Sigma_1, \Sigma_2$ be Pic-independent smooth projective curves. Then, for $q <c_2$ we have
$$H_q(\mathfrak{M}_{c_2}(\Sigma_1\times\Sigma_2),\mathfrak{M}^s_{c_2}(\Sigma_1\times\Sigma_2))=0.$$
\end{theorem}

Since we proved this result in the particular case of surfaces of general type which are product of curves, a natural question would be the following. 
\begin{question}
Does a  generalization of Thm.\thinspace\ref{lim}   holds for any surface of general type?
\end{question}

The question can also be stated in further generality for any surface  and furthermore for varieties of higher dimension. 
The underlying theme is a comparison of Kuranishi deformation theory as it described the local structure of moduli spaces and  moduli stacks. 

A more fundamental question is to actually determine the geometry of moduli spaces of vector bundles or instantons.  \cite{BEG} 
discusses the fundamental group of the moduli spaces of instantons on rational surfaces, 
but, a general result in this direction is lacking. So,
we end up
this note asking about the geometry of moduli on a complex surface $X$.

\begin{question}
What is the fundamental group of the moduli space of instantons (or holomorphic bundles) on $X$?
\end{question}

Certainly, we would like to know  higher homotopy groups as well, but it is important to clarify that even the 
knowledge of the fundamental groups of instanton moduli is still missing. 
 
 Finally, on a related topic, note that the local version of a moduli space is a deformation space. 
The reader interested on aspects of deformation theory might also consider reading the list of  20 open questions 
we mentioned in \cite{BGR}.
 
\paragraph{\bf Acknowledgements} 
 
 E. B. was partially supported by MIUR and GNSAGA of INdAM (Italy). E.G. and F.R. were supported by  the Vicerrector\'ia  de Investigaci\'on y Desarrollo Tecnol\'ogico UCN (Chile). F. R. was supported by Beca Doctorado Nacional
 Conicyt Folio 21170589 (Chile).

\end{document}